\documentclass[10pt]{article}
\usepackage{csit-2018}

\begin{document}
\English
\title
    {Non-Archimedean dynamics of the complex shift}
\author
    {L.\,B.~Tyapaev}
\email
    {tiapaevlb@info.sgu.ru}
\organization
    {Chernyshevsky Saratov State University, Saratov, Russia}
% если автор у статьи один, то команда \authorrefmark{1} не нужна.

\abstract
    {An (asynchronous) automaton transformation of one-sided infinite words over $p$-letter alphabet $\mathbb F_p=\mathbb Z/p\mathbb Z$, where $p$ is a prime, is a continuous transformation (w.r.t. the $p$-adic metric) of the ring of $p$-adic integers $\mathbb Z_p$. Moreover, an automaton mapping generates a non-Archimedean dynamical system on $\mathbb Z_p$. Measure-preservation and ergodicity (w.r.t. the Haar measure) of such dynamical systems play an important role in cryptography (e.g., in stream cyphers). The aim of this paper is to present a novel way of realizing a complex shift in $p$-adics. In particular, we introduce conditions on the Mahler expansion of a transformation on the $p$-adics which are sufficient for it to be complex shift. Moreover, we have a sufficient condition of ergodicity of such mappings in terms of Mahler expansion.}

\maketitle

\section*{Introduction}
Let $X$ be a finite set; we call this set an  {\it alphabet}. Given alphabet $X$, we denote by $X^*$ a free monoid generated by the set $X$. The elements of the monoid $X^*$ are expressed as words $x_0x_1\ldots x_{n-1}$ (including the empty word $\varnothing$). If $u=x_0x_1\ldots x_{n-1}\in X^*$, then $|u|=n$ is the length of the word $u$. The length of  $\varnothing$ is equal to zero. Along with finite words from $X^*$ we also consider infinite words of the form $x_0x_1x_2\ldots$, where $x_i\in X$. The set of such infinite words is denote by $X^{\infty}$. For arbitrary $u\in X^*$ and $v\in X^*\cup X^{\infty}$, we naturally defines the product (concatenation) $uv\in X^{\infty}$. A word $u\in X^*$ is the {\it beginning}, or {\it prefix} of a word $w\in X^{*}$ ($\in X^{\infty}$) if $w=uv$ for a certain $v\in X^*$ ($\in X^{\infty}$). The set $X^{\infty}$ is an infinite Cartesian product $X^{\mathbb N}$. We can introduce on the $X^{\mathbb N}$ the topology of the direct Tikhonov product of finite discrete topological spaces $X$. In this topology $X^{\infty}$ is homeomorphic to the Cantor set. Given finite word $u\in\ X^*$, the set $uX^{\infty}$ of all words beginning with $u$ is closed and open simultaneously (i.e. is {\it clopen}) in the given topology; the family of all such sets $\{uX^{\infty}: u\in X^*\}$ is the base of the topology.

We put a metric $d_{\pi}$ on $X^{\infty}$ by fixing a number $\pi>1$ and setting
$d_{\pi}(u,v)=\pi^{-\ell}$, where $\ell$ is the length of the longest common prefix of the words $u$ and $v$. The distance between identical words is equal to zero.

Thus the more that the initial terms of $u$ and $v$ agree, the closer they are to one another. It is easy to check that $d_{\pi}$ is a metric, and indeed a {\it non-Archimedean metric}, i.e. for any $u,v,w\in X^{\infty}$:
$$0\le d_{\pi}(u,v)\le 1;$$
$$d_{\pi}(u,v)=0 \Leftrightarrow u=v;$$
$$d_{\pi}(u,w)\le\max{\{d_{\pi}(u,v),d_{\pi}(v,w)\}}.$$
The set $uX^{\infty}$ of infinite words beginning with $u$ is a ball $B_{\pi^{-|u|}}(w)$ of radius $\pi^{-|u|}$ centered at arbitrary $w\in uX^{\infty}$.

Speaking about an {\it asynchronous automaton} we always understand the ``letter-to-word'' transducer $\mathfrak A=(\mathbb X,\mathcal S,\mathbb Y, h,g, s_0)$, where
\begin{enumerate}
\item[1)]$\mathbb X$ and $\mathbb Y$ are finite sets (the {\it input and output alphabets}, respectively);
\item[2)] $\mathcal S$ is a set ({\it the set of internal states of automaton});
\item[3)] $h\colon \mathbb X\times\mathcal S\to\mathcal S$ is a mapping ({\it transition function});
\item[4)] $g\colon\mathbb X\times\mathcal S\to\mathbb Y^*$ is a mapping ({\it output function}), and
\item[5)] $s_0\in\mathcal S$ is fixed ({\it initial state}).
\end{enumerate}

We assume that an asynchronous automaton $\mathfrak A$ works in a framework of discrete time steps. The automaton reads one symbol at a time, changing its internal state and outputting a finite sequence of symbols at each step.

The cardinality $\#\mathcal S$ of the set of states $\mathcal S$ of an automaton $\mathfrak A$ is called the cardinality of the automaton. In particular, automaton $\mathfrak A$ is finite if $\#\mathcal S<\infty$. If every value of the output function $g(\cdot,\cdot)$ is a one-letter word, then automaton $\mathfrak A$ is called a {\it synchronous} automaton. In the sequel denote a synchronous automaton via $\mathfrak B$. The functions $h$ and $g$ can be continued to the set $\mathbb X^*\times\mathcal S$ ($\mathbb X^{\infty}\times \mathcal S$).

The state $s\in\mathcal S$ of the automaton $\mathfrak A$ is called {\it accessible} if there exists a word $w\in \mathbb X^*$ such that $h(w,s_0)=s$. An automaton is called {\it accessible} if all its states are accessible. In the sequel, we consider only accessible automata.

An asynchronous automaton $\mathfrak A$ is {\it nondegenerate} if and only if there do not exist any accessible state $s\in\mathcal S$ and an infinite word $u\in \mathbb X^{\infty}$ such that, for an arbitrary prefix $w$ of the word $u$, the word $g(w,s)$ is empty.

A mapping $f\colon \mathbb X^{\infty}\to \mathbb Y^{\infty}$ is said to be defined a nondegenerate automaton $\mathfrak A$ if $f(u)=g(u,s_0)\in \mathbb Y^{\infty}$ for any $u\in \mathbb X^{\infty}$. {\it The mapping $f\colon\mathbb X^{\infty}\to\mathbb Y^{\infty}$ is continuous if and only if it is defined by a certain nondegenerate asynchronous automaton}, see \cite{Grigorchuk}. In the general case an asynchronous automaton defining a continuous mapping is infinite.

\section{$p$-adic numbers}
Let $p$ be a fixed prime number. By the fundamental theorem of arithmetics, each non-zero integer $n$ can be written uniquely as
$$n=p^{{\rm ord}_p n}\hat n,$$
where $\hat n$ is a non-zero integer, $p\nmid\hat n$, and ${\rm ord}_p n$ is a unique non-negative integer. The function ${\rm ord}_p n\colon \mathbb Z\backslash \{0\}\to\mathbb N_0$ is called the {\it $p$-adic valuation}. If $n,m\in \mathbb Z$, $m\ne0$, then the $p$-adic valuation of $x=n/m\in \mathbb Q$ as 
$${\rm ord}_p x={\rm ord}_p n-{\rm ord}_p m.$$
The $p$-adic valuation on $\mathbb Q$ is well defined; i.e., that ${\rm ord}_p x$ of $x$ does not depend on the fractional representation of $x$.

By using the $p$-adic valuation we will define a new absolute value on the
field $\mathbb Q$ of rational numbers. The $p$-adic absolute value of $x\in\mathbb Q\backslash \{0\}$ is given by
$$|x|_p=p^{-{\rm ord}_p x}$$
and $|0|_p=0$.
The $p$-adic absolute value is non-Archimedean. It induces the {\it $p$-adic
metric}
$$d_p(x,y)=|x-y|_p$$
which is non-Archimedean.
The completion of $\mathbb Q$ w.r.t. $p$-adic metric is a field, the {\it field of $p$-adic numbers, $\mathbb Q_p$.} The $p$-adic absolute value is extended to $\mathbb Q_p$, and $\mathbb Q$ is dense in $\mathbb Q_p$. The space $\mathbb Q_p$ is locally compact as topological space.

As the absolute value $|\cdot|_p$ may be only $p^{k}$ for some $k\in \mathbb Z$, for $p$-adic balls (i.e., for balls in $\mathbb Q_p$) we see that
$$B^{-}_{p^{k}}(x)=\{z\in\mathbb Q_p: d_p(z,x)<p^{k}\}=B_{p^{k-1}}(x)=\{z\in\mathbb Q_p: d_p(z,x)\le p^{k-1}\}.$$
Thus, we conclude that $p$-adic balls (of non-zero radii) are open and closed simultaneously; so $B_{p^k}(x)$ is a clopen ball of radius $p^k$ centered at $x\in \mathbb Q_p$.
Balls are compact; the set of all balls (of non-zero radii) form a topological base of a topology of a metric space. Thus, $\mathbb Q_p$ is a totally disconnected topological space.

The ball $B_1(0)=\{x\in\mathbb Q_p:|x|_p\le 1\}$ is {\it the ring of $p$-adic integers} and denoted via $\mathbb Z_p$. The space $\mathbb Z_p$ is a compact clopen totally disconnected metric subspace of $\mathbb Q_p$.

The ball $B^{-}_1(0)=\{x\in\mathbb Z_p:|x|_p<1\}=B_{p^{-1}}(0)=p\mathbb Z_p$ is a maximal ideal of the ring $\mathbb Z_p$. The factor ring 
$$\mathbb Z_p/B^{-}_1(0)=\mathbb Z/p\mathbb Z=\mathbb F_p$$ 
is then a finite field $\mathbb F_p$ of $p$ elements; it is called {\it the residue (class) field} of $\mathbb Q_p$. 

Every $p$-adic number has a unique representation as a sum of a special convergent $p$-adic series which is called a {\it canonical representation}, or a {\it $p$-adic expansion.} For $x\in \mathbb Z_p$ there exist a unique sequence $\delta_0(x),\delta_1(x),\ldots\in \{0,1,\ldots,p-1\}$ such that $x=\sum_{i=0}^{\infty}\delta_i(x)\cdot p^i=\delta_0(x)+\delta_1(x)\cdot p+\delta_2(x)\cdot p^2+\ldots$. The function $\delta_i(x)$ is called the {\it $i$-th coordinate function}.

Let $x=\delta_0(x)+\delta_1(x)\cdot p+\delta_2(x)\cdot p^2+\ldots$ be a $p$-adic integer in its canonical representation. The map
$$\bmod p^k: x=\sum_{i=0}^{\infty}\delta_i(x)\cdot p^i \mapsto x\bmod p^k=\sum_{i=0}^{k-1}\delta_i(x)\cdot p^i$$
is a {\it continuous ring epimorphism} of the ring $\mathbb Z_p$ onto the ring $\mathbb Z/p^k\mathbb Z$ of residues modulo $p^k$; it is called {\it reduction map modulo $p^k$}. The kernel of the epimorphism $\bmod p^k$ is a ball $p^k\mathbb Z_p=B_{p^{-k}}(0)$ of radius $p^{-k}$ around $0$. The rest $p^k-1$ balls of radii $p^{-k}$ are co-sets with respect to this epimorphism, e.g., $B_{p^{-k}}(1)=1+p^k\mathbb Z_p$, is a co-set of $1$, i.e. the sets of all $p$-adic integers
 that are congruent to $1$ modulo $p^k$.
 
 A $p$-adic integer $x\in\mathbb Z_p$ is {\it invertible} in $\mathbb Z_p$, that is, has a {\it multiplicative inverse} $x^{-1}\in \mathbb Z_p$, $x\cdot x^{-1}=1$, if and only if $\delta_0(x)\ne 0$; that is, if and only if $x$ is {\it invertible modulo $p$}, meaning $x\bmod p$ is invertible in $\mathbb F_p$. Invertible $p$-adic integer is also called a {\it unit}. The set $\mathbb Z_p^{*}$ of all units of $\mathbb Z_p$ is a group with respect to multiplication, called a {\it group of units}, or a {\it multiplicative subgroup} of $\mathbb Z_p$. The group of units $\mathbb Z_p^{*}$ is a $p$-adic {\it sphere} $S_1(0)$ of radius $1$ around $0$:
 $$\mathbb Z_p^{*}=\{z\in\mathbb Z_p:|z|_p=1\}=\mathbb Z_p\setminus p\mathbb Z_p=B_1(0)\setminus B_{p^{-1}}(0)=S_1(0).$$
 
Note that, the space $\mathbb Z_p$ is a {\it profinite algebra} with the structure of an {\it inverse limit}, that is the ring $\mathbb Z_p$ is an inverse limit of residue rings $\mathbb Z/p^{\ell}\mathbb Z$ modulo $p^{\ell}$, for all $\ell=1,2,3\ldots$
$$\mathbb Z_p\gets\ldots\gets\mathbb Z/p^2\mathbb Z\gets\mathbb Z/p\mathbb Z.$$

\section{Automata functions}

We identify $n$-letter words over $\mathbb F_p=\{0,1,\ldots,p-1\}$ with non-negative integers in a natural way: Given an $n$-letter word $u= x_0x_1\ldots x_{n-1}$, $x_i\in\mathbb F_p$, we consider $u$ as a base-$p$ expansion of the number $\alpha(u)=x_{n-1}\ldots x_1x_0=x_0+x_1\cdot p+\ldots+x_{n-1}\cdot p^{n-1}$. In turn, the latter number can be considered as an
element of the residue ring $\mathbb Z/p^{n}\mathbb Z=\{0,1,\ldots,p^n-1\}$ modulo $p^n$. Thus, every (synchronous) automaton $\mathfrak B=(\mathbb F_p,\mathcal S,\mathbb F_p,h,g,s_0)$ corresponds a map from $\mathbb Z/p^{n}\mathbb Z$ to $\mathbb Z/p^{n}\mathbb Z$, for every $n=1,2,3\ldots$.

Moreover, given an infinite word $u=x_0x_1x_2\ldots$ over $\mathbb F_p$ we consider $u$ as the $p$-adic integer $x=\alpha(u)=\ldots x_2x_1x_0$ whose canonical expansion is $x=\alpha(u)=x_0+x_1\cdot p+x_2\cdot p^2+\ldots=\sum_{i=0}^\infty{x_i\cdot p^i}$. Then every (synchronous) automaton $\mathfrak B=(\mathbb F_p,\mathcal S,\mathbb F_p,h,g,s_0)$ \textit{defines a map $f_{\mathfrak B}$ from the ring of $p$-adic integers $\mathbb Z_p$ to itself}: For every $x\in \mathbb Z_p$ we put $\delta_i(f_{\mathfrak B}(x))=g(\delta_i(x),s_i)$, $i=0,1,2,\ldots$ where $s_i=h(\delta_{i-1}(x),s_{i-1})$, $i=1,2,\ldots$. We say then that map $f_{\mathfrak B}$ is \textit{synchronous automaton function} (or, {\it automaton map}) of the synchronous automaton $\mathfrak B$.

Similarly way, a nondegenerate asynchronous automaton $\mathfrak A=(\mathbb F_p,\mathcal S,\mathbb F_p,h,g,s_0)$ naturally defines a {\it continuous mapping (w.r.t. $p$-adic metric)} $f_{\mathfrak A}$ from $\mathbb Z_p$ to $\mathbb Z_p$.

Note that, the automaton map $f_{\mathfrak B}\colon\mathbb Z_p\to\mathbb Z_p$ of the (synchronous) automaton $\mathfrak B=(\mathbb F_p,\mathcal S,\mathbb F_p,h,g,s_0)$ {\it satisfy the $p$-adic Lipschitz condition with constant } 1 (is a $1$-Lipschitz, for brevity), i.e. 
$$|f_{\mathfrak B}(x)-f_{\mathfrak B}(y)|_p\le|x-y|_p$$ for all $x,y\in\mathbb Z_p$. Conversely, {\it for every $1$-Lipschitz function $f\colon\mathbb Z_p\to\mathbb Z_p$ there exists an automaton $\mathfrak B=(\mathbb F_p,\mathcal S,\mathbb F_p,h,g,s_0)$ such that $f=f_{\mathfrak B}$}, see \cite {AnKh}. For example, every function defined by polynomial with $p$-adic integers coefficients (in particular, with rational integers) is a $1$-Lipschitz map, hence it is an automaton map.

Let $n\in \mathbb N$ be a natural number and let $\mathfrak C^{(n)}=(\mathbb X,\mathcal S,\mathbb Y,h,g,s_0)$ be a nongenerate asynchronous automaton, which is translated the infinite input word $u=x_0x_1\ldots x_{n-1}\ldots$  into infinite output word $w=\underbrace{\varnothing\varnothing\ldots\varnothing}_{n\text{~times} }y_ny_{n+1}\ldots$; So, we have 
$$y_i=g(x_i,s_i)=\varnothing \mbox{~for~}i=0,1,2\ldots,n-1,$$
$$s_i=h(x_{i-1},s_{i-1}) \mbox{~for~} i=1,2,\ldots,n-1, \mbox{~and}$$
$$y_i=g(x_i,s_i), s_{i+1}=h(x_i,s_i)$$ for all $i=n,n+1,\ldots$.

A {\it unilateral shift} is the transformation of the space of infinite words over alphabet $\mathbb X$ defined by the rule
$$x_0x_1x_2\ldots{\mapsto}~x_1x_2x_3\ldots.$$
Note that, a unilateral shift is defined by an asynchronous automaton $\mathfrak C^{(1)}=(\mathbb X,\mathcal S,\mathbb X,h,g,s_0)$, whose output function  $g$ is expressed as follows
$$g(x_0,s_0)=\varnothing, \mbox{~and~} g(x_i,s_i)=x_i \mbox{~for~$i=1,2,\ldots$}.$$
The initial state of this automaton, irrespective of the incoming letter outputs an empty word, after that the automaton outputs the incoming words without changes. In this case, the automaton $\mathfrak C^{(1)}=(\mathbb F_p,\mathcal S,\mathbb F_p,h,g,s_0)$ naturally defines the $p$-adic shift (or, {\it the one-sided Bernoulli shift}, see, e.g., \cite{Shift}); that is, the $p$-adic shift $\sigma\colon\mathbb Z_p\to\mathbb Z_p$ is expressed as follows. If $x=\ldots x_2x_1x_0=x_0+x_1p+x_2p^2+\ldots$, where $x_i\in\mathbb F_p$, we let $\sigma(x)=\frac{x-x_0}{p}=\ldots x_3x_2x_1=x_1+x_2p+x_3p^2+\cdots$. In other words, the shift $\sigma$ cuts off the first digit term in the $p$-adic expansion of $x\in\mathbb Z_p$.We see that if $\sigma^n$ denotes the $n$-fold iterate of $\sigma$, then we have $\sigma^n(x)=\frac{x-(x_0+x_1p+\ldots+x_{n-1}p^{n-1})}{p^n}=x_n+x_{n+1}p+\ldots$. Moreover, for $x\in \mathbb Z$, it is the case that $\sigma^n(x)=\lfloor\frac{x}{p^n}\rfloor$ were $\lfloor\cdot\rfloor$ is the greatest integer function.

A function $T\colon \mathbb Z_p\to \mathbb Q_p$ is called a {\it locally constant} if for every $x\in\mathbb Z_p$ there exist an open neighbourhood $U_x$ (e.g., a ball of radius $p^{-N}$ for some $N\in\mathbb N$ centered at $x$, $U_x=\{z\in\mathbb Z_p:|x-z|_p<p^{-N}\}$) such that $T$ is a constant on $U_x$. For example, for arbitrary $i\in\mathbb N$ a function $\delta_i(x)$ is locally constant, because $\delta_i$ remains unchanged if we replace $x$ by any $y$, such that $|x-y|_p<p^{-i}$.

Let $D\subset\mathbb Z_p$, not necessarily compact. A function $T\colon\mathbb Z_p\to\mathbb Q_p$ is called a {\it step function} on $D$ if there exists a positive integer $\ell$ such that $T(x)=T(y)$ for all  $x,y\in D$ with $|x-y|_p\le p^{-\ell}$. The smallest integer $\ell$ with this property is called the {\it order} of the step function $T$.

It is clear from the definition that a step function is a locally constant on $D$. On $\mathbb Z_p$ it also holds, that any locally constant function is a step function.

Let $f\colon\mathbb Z_p\to\mathbb Z_p$ be a 1-Lipschitz function. Given natural $n\in \mathbb N$, for all $x\in\mathbb Z_p$ we can represent $f$ as
$$f(x)=(f(x\bmod p^n))\bmod p^n+p^nG_z(t),$$
where $t=p^{-n}(x-(x\bmod p^n))\in\mathbb Z_p$, $z=x\bmod p^n$, and $G_z\colon \mathbb Z_p\to\mathbb Z_p$ is a 1-Lipschitz function. It is clear that $f$ is the automaton map of a synchronous automaton $\mathfrak B=(\mathbb F_p,\mathcal S,\mathbb F_p,h,g,s_0)$, and $G_z$ is an automaton map of the automaton $\mathfrak B_z=(\mathbb F_p,\mathcal S,\mathbb F_p,h,g,s(z))$, where $s(z)\in\mathcal S$ is the accessible state of the automaton $\mathfrak B$, i.e. $s(z)$ was reached after being feeded by the input word $z=x\bmod p^n$.

Similarly, for map $f\colon\mathbb Z_p\to\mathbb Z_p$ that defined by an asynchronous automaton $\mathfrak C^{(n)}$, we can see that
$$f(x)=G_{z}(t)+T(x),$$
where for  any $z=x\bmod p^n$, the map $G_{z}\colon\mathbb Z_p\to\mathbb Z_p$ is a 1-Lipschitz, $T(x)$ is a step function of order not greater than $n$, and $t=p^{-n}(x-(x\bmod p^n))$; and we say that $f$ is a {\it complex shift}.

For a complex shift $f$ a following condition holds: There exist positive integer $M\in \mathbb N$ such that for every $i\ge M$ the $i$-th coordinate function $\delta_i(f(x))$ does not depend on $\delta_{i+k}(x)$ for $k=1,2,\ldots$. Hence, a complex shift is a {\it locally 1-Lipschitz function}. By the definition, a function $F$ is a locally $1$-Lipschitz if  for a given $x\in\mathbb Z_p$, there exist an open neighbourhood $U_x$ of $x$ such that the inequality
$$|F(x)-F(y)|_p\le|x-y|_p$$
holds for all $y\in U_x$. As $\mathbb Z_p$ is compact, the function $F\colon\mathbb Z_p\to\mathbb Z_p$ is a locally $1$-Lipschitz if and only if the latter inequality holds for all $x,y\in\mathbb Z_p$ which are sufficiently close to one another.

\section{Dynamical systems}
A \textit {dynamical system} on a measurable space $\mathbb S$ is understood as a triple $(\mathbb S,\mu,f)$, where $\mathbb S$ is a set endowed with a measure $\mu$ and $f\colon\mathbb S\to \mathbb S$ is a measurable function;  that is, the $f$-preimage $f^{-1}(T)$ of any $\mu$-measurable subset $T\subset \mathbb S$ is a $\mu$-measurable subset of $\mathbb S$.

An iteration of a function $f_{\mathfrak A}\colon\mathbb Z_p\to\mathbb Z_p$  which is defined by (asynchronous) automaton $\mathfrak A=(\mathbb F_p,\mathcal S,\mathbb F_p,h,g,s_0)$ generates a dynamical system $(\mathbb Z_p,\mu_p, f_{\mathfrak A})$ on the space $\mathbb Z_p$. The space $\mathbb Z_p$ is equipped with a natural probability measure, namely, the {\it Haar measure} $\mu_p$ normalized so that the measure of the whole space is 1, $\mu_p(\mathbb Z_p)=1$. Balls $B_{{p}^{-k}}(a)$ of nonzero radii constitute the base of the corresponding $\sigma$-algebra of measurable subsets of $\mathbb Z_p$. That is, every element of the $\sigma$-algebra, the measurable subset of $\mathbb Z_p$, can be constructed from the elementary measurable subsets by taking complements and countable unions. We put $\mu_p(B_{{p}^{-k}}(a))=p^{-k}$. 

We remind that if a measure space $\mathbb S$ endowed with a probability measure $\mu$ is also a topological space, the measure $\mu$ is called {\it Borel} if all Borel sets in $\mathbb S$ are $\mu$-measurable. Recall that a {\it Borel set} is any element of $\sigma$-algebra generated by all open subsets of $\mathbb S$; that is, a Borel subset can be constructed from open subsets with the use of complements and countable unions. A probability measure $\mu$ is called {\it regular} if for all Borel sets $X$ in $\mathbb S$
$$\mu(X)=\sup\{\mu(A): A\subseteq X, A \mbox{~closed}\}=\inf\{\mu(B): X\subseteq B, B\mbox{~open}\}.$$
The probability measure $\mu_p$ is Borel and regular.

A dynamical system $(\mathbb Z_p,\mu_p, f_{\mathfrak A})$ is also {\it topological} since $\mathbb Z_p$ are not only measurable space but also metric space, and corresponding transformation $f_{\mathfrak A}$ are not only measurable but also continuous. Moreover, this dynamical system is non-Archimedean, due to the fact that the space $\mathbb Z_p$ is non-Archimedean space.

A measurable mapping $f\colon\mathbb Z_p\to \mathbb Z_p$ is called {\it measure-preserving} if $\mu_p(f^{-1}(S))=\mu_p(S)$ for each measurable subset $S\subset\mathbb Z_p$.
A measure-preserving map $f$ is said to be \textit{ergodic} if for
each measurable subset $S$ such that $f^{-1}(S)=S$ holds either $\mu_p(S)=1$ or $\mu_p(S)=0$; so ergodicity of the map $f$ just means that $f$ has no proper invariant subsets; that is, invariant subsets whose measure is neither $0$ nor $1$.

The following question arises. What continuous with respect to the metric $\mu_p$ transformations are measure-preserving or ergodic with respect to the mentioned measure?

For a given $f\colon\mathbb Z_p\to\mathbb Z_p$ and $n\in\mathbb N$, let $f_k$ be a function defined on the ring $\mathbb Z/p^{n\cdot k}\mathbb Z$ and valuated in the ring $\mathbb Z/p^{n\cdot {(k-1)}}\mathbb Z$, where $k=2,3,\ldots$.

The following criterion of measure-preservation for a complex shift $f$ is valid: {\it A mapping $f\colon\mathbb Z_p \to \mathbb Z_p$ is measure-preserving if and only if the number $\#f_k^{-1}(x)$ of $f_k$-preimages of the point $x\in \mathbb Z/p^{n\cdot(k-1)}\mathbb Z$ is equal to $p^n$} \cite {T11,T12}.

Given a map $f\colon\mathbb Z_p\to\mathbb Z_p$, a point $z_0\in \mathbb Z_p$ is said to be a \textit{periodic point} if there exists $r\in\mathbb N$ such that $f^r(z_0)=z_0$. The least $r$ with this property is called the \textit{length} of period of $z_0$. If $z_0$ has period $r$, it is called an \textit{$r$-periodic point}. The orbit of an $r$-periodic point $z_0$ is $\{z_0, f(z_0),\ldots,f^{r-1}(z_0)\}$. This orbit is called an \textit{$r$-cycle}.

For a given $n\in\mathbb N$, let $f\bmod p^{k\cdot n}\colon\mathbb Z/p^{k\cdot n}\mathbb Z\to\mathbb Z/p^{k\cdot n}\mathbb Z$, for $k=1,2,3,\ldots$; and  let $\gamma_k$ be an $r_k$-cycle $\{z_0,z_1,\ldots,z_{r_k-1}\}$, where $z_j=(f\bmod  p^{k\cdot n})^j(z_0)$, $0\le j\le r_k-1$, $k=1,2,3,\ldots$.

The following condition of ergodicity holds: {\it Let $f\colon \mathbb Z_p\to\mathbb Z_p$ be a complex shift and let $f$ be a measure-preserving map. Then $f$ is ergodic if for every $k\in \mathbb N$ $\gamma_k$ is a unique cycle} \cite {T13}.

\section{Measure-preservation and ergodicity in terms of Mahler expansion}

By Mahler's Theorem, any continuous function $F\colon \mathbb Z_p\to \mathbb Z_p$ can be expressed in the form of a uniformly convergent series, called its \emph {Mahler Expansion} (or, {\it Mahler series}):
$$F(x)=\sum_{m=0}^{\infty}a_m\binom{x}{m},$$
where $$a_m=\sum_{i=0}^{m}(-1)^{m+i}F(i)\binom{m}{i}\in\mathbb Z_p$$ and
$$\binom{x}{m}=\frac{x(x-m)\cdots(x-m+1)}{m!}, m=1,2,\ldots, \binom{x}{0}=1.$$
Mahler series converges uniformly on $\mathbb Z_p$ if and only if
$$\stackrel{p}{\lim_{m\to\infty}}a_m=0.$$
Hence uniformly convergent series defines a uniformly continuous function on $\mathbb Z_p$. The function $f$ represented by the Mahler series is uniformly differentiable everywhere on $\mathbb Z_p$ if and only if
$$\stackrel{p}{\lim_{m\to\infty}}{\frac{a_{m+k}}{m}}=0$$
for all $k\in\mathbb N$.

The function $f$ is analytic on $\mathbb Z_p$ if and only if
$$\stackrel{p}{\lim_{m\to\infty}}{\frac{a_m}{m!}}=0.$$

Various properties of the function $f$ can be expressed via
properties of coefficients of its Mahler expansion. 

Let $a_m^{(n)}$ be the $n$-th Mahler coefficient of the Bernoulli shift $\sigma^n$. We have
$$\sigma^n(x)=\sum_{m=0}^{\infty}a_m^{(n)}\binom{x}{m}.$$
The coefficients $a_m^{(n)}$ satisfy the following properties \cite{Shift}:{\it
\begin{enumerate}
\item []$a_m^{(n)}=0$ for $0\le m<p^n$;
\item []$a_m^{(n)}=1$ for $m=p^n$;
\item []Suppose $j\ge 0$. Then, $p^j$ divides $a_m^{(n)}$ for $m>jp^n-j+1$ (and so, $|a_m^{(n)}|_p\le 1/p^j$).
\end{enumerate}
}
The following statement gives a description of $1$-Lipschitz measure-preserving (respectively, of $1$-Lipschitz ergodic) transformations on $\mathbb Z_p$ \cite{AnKh}.
{\it The function $f$ defines a $1$-Lipschitz measure-preserving transformation on $\mathbb Z_p$ whenever the following conditions hold simultaneously:
\begin {itemize}
\item[] $a_1\not \equiv 0\pmod p$;
\item[] $a_m\equiv 0 \pmod {p^{\lfloor\log_p m\rfloor+1}}$, $m=2,3,\ldots$
\end{itemize}
The function $f$ defines a $1$-Lipschitz ergodic transformation on $\mathbb Z_p$ whenever the following conditions hold simultaneously:
\begin {itemize}
\item[] $a_0\not \equiv 0\pmod p$;
\item[] $a_1\equiv 1\pmod p$ for $p$ odd;
\item[] $a_1\equiv 1 \pmod 4$ for $p=2$;
\item[] $a_m\equiv 0 \pmod {p^{\lfloor\log_p (m+1)\rfloor+1}}$, $m=1,2,3,\ldots$
\end{itemize}
Moreover, in the case $p=2$ these conditions are necessary.}

The following statement gives a description of complex shift in terms of Mahler expansion.

\begin{Theorem}
 A function $f\colon \mathbb Z_p\to \mathbb Z_p$ is a complex shift if and only if
$$|a_m|_p\le p^{-\lfloor \log_{p^n}m\rfloor+1},$$ where $n\in \mathbb N$, $m\ge 1$.
\end{Theorem}

The following theorems give a description of measure-preservation and ergodicity for a complex shift.

\begin{Theorem}
A complex shift $f\colon\mathbb Z_p\to\mathbb Z_p$ is measure-preserving whenever the following conditions hold  simultaneously:
\begin{enumerate}
\item []$a_m\not \equiv 0\pmod p$ for $m=p^n$;
\item []$a_m\equiv 0 \pmod {p^{\lfloor \log_{p^n}m\rfloor}}$, $m>p^n$,
\end{enumerate}
where $n\in\mathbb N$.
\end{Theorem}

\begin{Theorem}
A complex shift $f\colon\mathbb Z_p\to\mathbb Z_p$ is ergodic on $\mathbb Z_p$ whenever the following conditions hold  simultaneously:
\begin{enumerate}
\item []$a_1+a_2+\ldots+a_{p^n-1}\equiv 0 \pmod p$;
\item []$a_m\equiv 1\pmod p$ for $m=p^n$;
\item []$a_m\equiv 0 \pmod {p^{\lfloor \log_{p^n}m\rfloor}}$, $m>p^n$,
\end{enumerate}
where $n\in\mathbb N$.
\end{Theorem}

Let $f\colon\mathbb Z_p\to\mathbb Z_p$ be a complex shift, and let $E_k(f)$ be a set
of all the following points $e_k^f (x)$ of Euclidean unit square $\mathbb I^2=[0, 1]\times[0, 1]\subset
\mathbb R^2$ for $k=1,2,3,\ldots$ \cite{arxiv}:
$$e_k^f (x)=\Bigl(\frac{x\bmod p^{n+k}}{p^{n+k}},\frac{f(x)\bmod p^{k}}{p^{k}}\Bigr),$$
where $x\in\mathbb Z_p$, $n\in\mathbb N$.
Note that $x \bmod p^{n+k}$ corresponds to the prefix of length $n+k$ of the infinite word $x\in\mathbb Z_p$, i.e., to the input word of length $n+k$ of the automaton $\mathfrak C^{(n)}$; while $f(x) \bmod p^k$ corresponds to the respective output word of length $k$. Denote via $\mathcal E(f)$ the closure of the set $E(f)=\bigcup_{k=1}^{\infty}E_k(f)$ in the topology of real plane $\mathbb R^2$. As $\mathcal E(f)$ is closed, it is measurable with respect to the Lebesgue measure on real plane $\mathbb R^2$. Let $\lambda(f)$ be the Lebesgue measure of $\mathcal E(f)$.

\begin{Theorem}
For a given complex shift $f\colon\mathbb Z_p\to\mathbb Z_p$ the closure $\mathcal E(f)$ is nowhere dense in $\mathbb I^2$, hence $\lambda(f)=0$.
\end{Theorem}

\end{document}